\documentclass{amsart}
\usepackage{}
\usepackage{amsmath}
\usepackage{amscd}
\usepackage{amssymb}
\usepackage{amsfonts}
\usepackage{hyperref}
\newtheorem{theorem}{Theorem}[section]
\newtheorem{lemma}[theorem]{Lemma}

\theoremstyle{definition}

\theoremstyle{remark}
\newtheorem{remark}[theorem]{Remark}
\numberwithin{equation}{section}
\begin{document}

\title[Sharp Hamilton's Laplacian estimate for the heat kernel]
{Sharp Hamilton's Laplacian estimate for the heat kernel
on complete manifolds}

\author{Jia-Yong Wu}
\address{Department of Mathematics, Shanghai Maritime University,
Haigang Avenue 1550, Shanghai 201306, P. R. China}

\email{jywu81@yahoo.com}

\thanks{This work is partially supported by the NSFC (No.11101267).}

\subjclass[2000]{Primary 58J35; Secondary 35K08.}
\date{\today}

\dedicatory{}

\keywords{gradient estimate, heat kernel,
heat equation.}

\begin{abstract}
In this paper we give Hamilton's Laplacian estimates for the heat
equation on complete noncompact manifolds with nonnegative Ricci
curvature. As an application, combining Li-Yau's lower and upper
bounds of the heat kernel, we give an estimate on Laplacian form
of the heat kernel on complete manifolds with nonnegative Ricci
curvature that is sharp in the order of time parameter for the heat
kernel on the Euclidean space.
\end{abstract}
\maketitle
\section{Introduction}\label{sec1}
In \cite{Hamilton}, R. Hamilton established the following gradient and
Laplacian estimates for a bounded positive solution to the heat equation
on closed manifolds with Ricci curvature bounded below.
\begin{theorem}[R. Hamilton \cite{Hamilton}]\label{ham1}
Let $(M,g)$ be an $n$-dimensional closed Riemannian manifold
with Ricci curvature satisfying $Ric\geq-K$ for some constant
$K\geq 0$. Assume that $u$ is a positive solution to the heat
equation with $u\leq A$ on $M\times[0,T]$ for some constant
$A<\infty$, where $0<T<\infty$. Then
\begin{equation}\label{gradest1}
t\frac{|\nabla u|^2}{u^2}\leq(1+2Kt)\log\left(\frac{A}{u}\right).
\end{equation}
If we further assume $T\leq1$, then for $0\leq t\leq T$,
\begin{equation}\label{gradest2}
t\frac{\Delta u}{u}
\leq C(n,K)\left[1+\log\left(\frac{A}{u}\right)\right].
\end{equation}
\end{theorem}

For estimate \eqref{gradest2}, when $K=0$, by choosing function
$\varphi=t$ in Hamilton's proof of Lemma 4.1 in \cite{Hamilton},
we easily confirm that the condition ``$T\leq1$" can be removed.
Hamilton's estimate \eqref{gradest1} shows that one can compare
two different points at the same time, however the well-known
Li-Yau's gradient estimate \cite{[Li-Yau]} only allows comparisons
between different points at the different times. In \cite{Kotsch},
B. Kotschwar generalized the gradient estimate \eqref{gradest1}
to the case of complete, noncompact Riemannian manifolds with
Ricci curvature bounded below. Using this generalization,
Kotschwar gave an estimate on the gradient of the heat
kernel for complete manifolds with nonnegative Ricci curvature.
Moreover this estimate is sharp in the order of $t$ for the
heat kernel on $\mathbb{R}^n$.

Kotschwar's result can be used to prove the monotonicity of Ni's
entropy functional (stationary metric for Perelman's
$\mathcal{W}$-functional in \cite{[Perelman]}) in \cite{[Ni1]}
for the fundamental solution to the heat equation on complete,
noncompact manifolds. We would like to point out that, in the
course of justifying the monotonicity Ni's entropy functional
on complete noncompact manifolds, one may need a noncompact
version of Hamilton's Laplacian estimate \eqref{gradest2}.
However, as far as we know, perhaps no people generalized
the estimate \eqref{gradest2} to the complete noncompact case.
In \cite{[CCG]}, the authors only briefly sketched a proof of
the estimate \eqref{gradest2} for the noncompact case but
didn't give any detail. In this paper, we will provide a
full detailed proof that Hamilton's Laplacian estimate
\eqref{gradest2} also holds for complete, noncompact
manifolds with nonnegative Ricci curvature. Precisely,
we show that
\begin{theorem}\label{the1}
Let $(M,g)$ be an $n$-dimensional complete noncompact Riemannian
manifold with nonnegative Ricci curvature. Suppose $u$ is a
smooth positive solution to the heat equation
\begin{equation}\label{weheat}
\frac{\partial u}{\partial t}-\Delta u=0,
\end{equation}
satisfying $u\leq A$ for some constant $A<\infty$ on
$M\times[0,T]$, where $0<T<\infty$.
Then
\begin{equation}\label{gradest4}
t\frac{\Delta u}{u}
\leq n+4\log\left(\frac{A}{u}\right)
\end{equation}
for all $x\in M$ and $0\leq t\leq T$.
\end{theorem}

The proof of Theorem \ref{the1} is similar to the arguments of
Kotschwar \cite{Kotsch}, which can be divided into two steps.
In the first step, we obtain some Bernstein-type estimate of
$\Delta u$, similar to the upper estimate of $|\nabla u|$
derived by Kotschwar \cite{Kotsch} on complete noncompact
manifolds. In the second step, using upper estimates of
$\Delta u$ and $|\nabla u|$, we apply the maximum principle
to the quantity of Hamilton's Laplacian estimate on complete
noncompact Riemannian manifolds due to Karp-Ni \cite{KaLi}
or Ni-Tam \cite{[Ni-Tam]}. We remark that a priori integral
bound needed for the application of maximum principle on
complete noncompact manifolds has also been obtained in
\cite{[Gri]} and in \cite{Yu} for more general setting.

As an application of Theorem \ref{the1}, we obtain the following
Laplacian estimate of the heat kernel on a complete noncompact
Riemannian manifold with nonnegative Ricci curvature.
\begin{theorem}\label{Kethe1}
Let $(M,g)$ be an $n$-dimensional complete noncompact Riemannian manifold
with nonnegative Ricci curvature, and $H(x,y,t)$ its heat kernel.
Then, for all $\delta>0$, there exists a constant $C=C(n,\delta)$
such that
\[
\frac{\Delta H}{H}(x,y,t)\leq\frac{2}{t}
\left[C+4\frac{d^2(x,y)}{(4-\delta)t}\right]
\]
for all $x,y\in M$ and $t>0$.
\end{theorem}
\begin{remark}
We would like to point out that Theorem \ref{Kethe1} is sharp
in the order of $t$ for the heat kernel on $\mathbb{R}^n$.
\end{remark}

The structure of this paper is organized as follows. In Section
\ref{sec2}, we derive Bernstein-type gradient estimates of the
Laplacian for solutions to the heat equation (see Theorem
\ref{theor2}). Our proof makes use of Shi's gradient estimates
\cite{[Shi]}, combining the classical cut-off function arguments.
In Section \ref{sect3}, we finish the proof of Theorem \ref{the1}
by using Theorem \ref{theor2}. In Section \ref{sect4}, we apply
Theorem \ref{the1} to the heat kernel and complete the proof
of Theorem \ref{Kethe1}.

\section{Bernstein-type estimates}\label{sec2}
In this section, we assume that $(M,g)$ be an $n$-dimensional
complete noncompact Riemannian manifold with the Ricci curvature
uniformly bounded below by $-K$ for some constant $K\geq 0$,
and suppose that $u$ is a smooth solution to the heat equation
\eqref{weheat} satisfying $|u|\leq A$ on some open $U\subset M$
for $0\leq t\leq T<\infty$. At first we recall the Kotschwar's
result in \cite{Kotsch}.
\begin{theorem}[Kotschwar \cite{Kotsch}]\label{Ber1}
Let $(M,g)$ be an $n$-dimensional complete noncompact Riemannian
manifold with $Ric\geq-K$ for some constant $K\geq 0$. Suppose
$u$ is a smooth solution to the heat equation \eqref{weheat}
satisfying $|u|\leq A$ on $B_p(2R)\times[0,T]$ for some
$p\in M^n$ and $A,R,T>0$. Then there exists a constant
$C=C(n,K)$ such that
\[
t|\nabla u|^2\leq CA^2\left[1+T\left(1+\frac{1}{R^2}\right)\right]
\]
holds on $B_p(R)\times[0,T]$.
\end{theorem}

\begin{remark}\label{rem1}
If $Ric\geq0$, from the proof course of Theorem \ref{Ber1}
in \cite{Kotsch}, one shows that
\[
t|\nabla u|^2\leq C(n)A^2\left(1+\frac{T}{R^2}\right)
\]
on $B_p(R)\times[0,T]$.
\end{remark}

\begin{remark}
Letting $R\to \infty$ in the proof course of Theorem \ref{Ber1},
one immediately shows that there exists a constant $C(n)$
such that
\begin{equation}\label{td1}
t|\nabla u|^2\leq C(n)A^2(1+KT)
\end{equation}
on $M^n\times[0,T]$.
\end{remark}

In the above description, Kotschwar showed the first derivative
estimate of the positive solution to the heat equation on
complete manifolds. Below we will give an upper estimate of
$\Delta u$. Our proof is similar in spirit to the derivative
estimates due to Shi \cite{[Shi]} (see also \cite{Kotsch}). Let
\begin{equation}\label{deffor}
F(x,t):=(C+t|\nabla u|^2)t^2|\Delta u|^2,
\end{equation}
where the constant $C$ is to be chosen. The following lemma
is useful for proving Theorem \ref{theor2}.
\begin{lemma}\label{lemq3}
Let $(M,g)$ be an $n$-dimensional complete Riemannian manifold with
$Ric\geq -K$ for some constant $K\geq 0$. If $0<u\leq A$ is the
solution to the heat equation \eqref{weheat} on $B_p(2R)\times[0,T]$
for some $p \in M$ and $A,R,T>0$, where $T\leq 1$, satisfying
\[
|\nabla u|^2\leq \frac{C_*}{t}
\]
for some constant $C_*$ on $B_p(R)\times(0,T]$, then
there exists a finite positive constant $c:=c(n,K,A,R)$ such that
\[
\frac{\partial F}{\partial t}
\leq \Delta F-\frac {c}{t}F^2+\frac{C_*^2}{t}
\]
on $B_p(R)\times(0,T]$.
\end{lemma}
\begin{remark}
The assumption $T\leq 1$ in Lemma \ref{lemq3} is only used
in \eqref{gjgj}. By Theorem \ref{Ber1}, one may choose
$C_*:=C(n,K)A^2\left[1+T\left(1+\frac{1}{R^2}\right)\right]$.
Also we can choose $c:=C^{-1}(n)\cdot C^{-2}_*$. Moreover,
if $R\to \infty$, then $\lim_{R\to \infty}c$ is a finite
positive constant. Note that here the constant $C(n,K)$
may be different from the one in Theorem \ref{Ber1}.
\end{remark}

\begin{remark}\label{rem26}
When $K=0$, from \eqref{gjgj}, we see that the assumption $T\leq 1$
can be replaced by $T<\infty$. In this case, from Remark
\ref{rem1}, one can choose $C_*:=C(n)A^2\left(1+\frac{T}{R^2}\right)$
and $c:=C^{-1}(n)\cdot C^{-2}_*$. If $R\to \infty$, then
$\lim_{R\to \infty}c$ still be a finite positive constant,
\textbf{independent} on $T$.
\end{remark}

\begin{proof}
At first, the evolution formula of $t|\nabla u|^2$ is that
\begin{equation*}
\begin{aligned}
\left(\frac{\partial}{\partial t}-\Delta \right)(t|\nabla u|^2)
&=-2t|\nabla\nabla u|^2-2t Ric(\nabla u,\nabla u)+|\nabla u|^2\\
&\leq-t|\nabla\nabla u|^2-\frac{t}{n}|\Delta u|^2-2t Ric(\nabla u,\nabla u)+|\nabla u|^2\\
&\leq-t|\nabla\nabla u|^2-\frac{t}{n}|\Delta u|^2+(2Kt+1)|\nabla u|^2,
\end{aligned}
\end{equation*}
where we used $Ric\geq -K$. Then we compute that
\[
\left(\frac{\partial}{\partial t}-\Delta\right)|\Delta u|^2
=-2|\nabla \Delta u|^2
\]
and hence
\[
\left(\frac{\partial}{\partial t}-\Delta\right)(t^2|\Delta u|^2)
=-2t^2|\nabla \Delta u|^2+2t|\Delta u|^2.
\]
By the assumption of this lemma, we can choose $C$ in \eqref{deffor}
such that $C=8C_*$, which implies that $8t|\nabla u|^2\leq C$.
Combining the above equations yields
\begin{equation*}
\begin{aligned}
\left(\frac{\partial}{\partial t}-\Delta\right)F
&=(C+t|\nabla u|^2)\left[\left(\frac{\partial}{\partial t}-\Delta\right)
(t^2|\Delta u|^2)\right]\\
&\quad+\left[\left(\frac{\partial}{\partial t}-\Delta\right)(C+t|\nabla u|^2)\right]t^2|\Delta u|^2\\
&\quad-2t^3\nabla(|\nabla u|^2)\cdot\nabla((\Delta u)^2)\\
&\leq(C+t|\nabla u|^2)\left(-2t^2|\nabla \Delta u|^2+2t|\Delta u|^2\right)\\
&\quad+\left[-t|\nabla\nabla u|^2-\frac{t}{n}|\Delta u|^2+(2Kt+1)|\nabla u|^2\right]t^2|\Delta u|^2\\
&\quad+8t^3|\nabla u||\nabla\nabla u|\cdot|\Delta u||\nabla \Delta u|\\
&\leq-18t^3|\nabla u|^2|\nabla \Delta u|^2-\frac{t^3}{n}(\Delta u)^4-t^3|\nabla\nabla u|^2(\Delta u)^2\\
&\quad+2t(\Delta u)^2(C+t|\nabla u|^2)+(2Kt+1)t^2|\nabla u|^2|\Delta u|^2\\
&\quad+t^3|\nabla\nabla u|^2\cdot|\Delta u|^2+16t^3|\nabla u|^2\cdot|\nabla \Delta u|^2,
\end{aligned}
\end{equation*}
where we used the Schwarz inequality. Since $t\leq T\leq1$, the above formula becomes
\begin{equation}
\begin{aligned}\label{gjgj}
\left(\frac{\partial}{\partial t}-\Delta\right)F
&\leq-\frac{t^3}{n}(\Delta u)^4+(2Kt+1)Ct|\Delta u|^2+4Ct(\Delta u)^2\\
&\leq-\frac{t^3}{n}(\Delta u)^4+(2K+1)Ct|\Delta u|^2+4Ct(\Delta u)^2\\
&\leq-\frac {c}{t}F^2+\frac{18n(1+K^2)C^2}{t},
\end{aligned}
\end{equation}
where in the last inequality, we used
\[
\frac{t^3}{2n}(\Delta u)^4+\frac{18n(1+K^2)C^2}{t}\geq C(2K+1)t|\Delta u|^2+4Ct(\Delta u)^2.
\]
Here $c:=c(n,K,A,R)$ depends on $n$, $K$, $A$ and $R$.
For example, we may choose $c:=C^{-1}(n)\cdot C^{-2}_*$.
Then the result follows.
\end{proof}

Using Lemma \ref{lemq3}, we prove the following Laplacian estimate
for the positive solution to the heat equation.
\begin{theorem}\label{theor2}
Let $(M,g)$ be an $n$-dimensional complete  Riemannian manifold with
$Ric\geq -K$ for some constant $K\geq 0$. Let $u$ be a positive
solution to the heat equation \eqref{weheat} with $u\leq A$
on $B_p(2R)\times(0,T]$ for some
$p\in M^n$ and $A,R,T>0$, where $A<\infty$ and $T\leq 1$.
Then there exists a constant $C=C(n,K)$ such that
\begin{equation}\label{td2}
t|\Delta u|\leq CA\left[1+T\left(1+\frac{1}{R^2}\right)\right]^{1/2}
\cdot\left(1+\frac{T}{R^2}\right)^{1/2}
\end{equation}
on $B_p(R)\times[0,T]$.
\end{theorem}

\begin{remark}\label{Laptd1}
When $K=0$, the assumption $T\leq 1$ can be replaced by
$T<\infty$. In this case, from Remark \ref{rem26},
estimate \eqref{td2} can be rewritten by a simple version
\begin{equation}\label{td2gen}
t|\Delta u|\leq C(n)A\left(1+\frac{T}{R^2}\right)
\end{equation}
on $B_p(R)\times[0,T]$. If we further let $R\to \infty$,
then
\begin{equation}\label{td3gen}
t|\Delta u|\leq C(n)A
\end{equation}
on $M\times[0,T]$.
\end{remark}

\begin{proof}[Proof of Theorem \ref{theor2}]
As in \cite{[Li-Yau],[Cala]} (see also \cite{Kotsch} or \cite{Wu10}),
for any $p\in M$ and $R>0$,
we may choose a cut-off function with $\eta(x)=1$ on $B_p(R)$
and supported in $B_p(2R)$ satisfying
\begin{equation}
|\nabla \eta|^2\leq \frac{C_3}{R^2}\eta
\end{equation}
and
\begin{equation}
\Delta \eta\geq-\frac{C_3}{R^2}
\end{equation}
for some $C_3=C_3(n)>0$. Letting $G:=\eta F$, we compute that
\[
\left(\frac{\partial}{\partial t}-\Delta\right)G \leq
\frac{\eta}{t}(-cF^2+C^2_*)-\Delta\eta\cdot F-2\nabla \eta\cdot\nabla F,
\]
where $C_*=C(n,K)A^2\left[1+T\left(1+\frac{1}{R^2}\right)\right]$.
Assume that at a point $(x_0,t_0)$ where the function $G$ attains
its positive maximum in $B(p,2R)\times(0,T]$. Then at $(x_0,t_0)$
we have
\[
\left(\frac{\partial}{\partial t}-\Delta\right)G\geq 0
\quad \mathrm{and}\quad 0=\nabla G=\eta \nabla F+F\nabla\eta.
\]
Therefore at $(x_0,t_0)$, we have
\begin{equation*}
\begin{aligned}
0&\leq\frac{1}{t}\left[-c(\eta F)^2+C^2_*\eta\right]-\Delta\eta\cdot (\eta F)+2F|\nabla\eta|^2\\
&\leq\frac{1}{t}[-cG^2+C^2_*\eta]+\frac{C_3}{R^2}G+2\frac{C_3}{R^2}G
\end{aligned}
\end{equation*}
and hence
\begin{equation}\label{guanji}
cG^2(x_0,t_0)\leq C^2_*+3\frac{C_3}{R^2}G(x_0,t_0)t_0.
\end{equation}
Since
\[
3\frac{C_3}{R^2}G(x_0,t_0)t_0\leq \frac c2G^2(x_0,t_0)+\frac{8C_3^2}{cR^4}t^2_0,
\]
inequality \eqref{guanji} implies that
\[
\frac c2G^2(x_0,t_0)\leq C^2_*+\frac{8C_3^2}{cR^4}t^2_0.
\]
Therefore for any $(x,t)\in B_p(R)\times[0,T]$,
\[
G^2(x,t)\leq G^2(x_0,t_0)\leq 2c^{-1}C^2_*+c^{-2}\frac{16C_3^2}{R^4}t^2_0.
\]
Since $c:=C^{-1}(n)\cdot C^{-2}_*$, we have that
\[
G^2(x,t)\leq C(n)C^4_*+\frac{C(n)}{R^4}T^2C^4_*
\]
for any $(x,t)\in B_p(R)\times[0,T]$. This implies
\[
G(x,t)\leq C(n)C^2_*+\frac{C(n)}{R^2}TC^2_*
\]
for any $(x,t)\in B_p(R)\times[0,T]$. By the definitions of
$C_*$ and $G$, we have
\[
8C_*t^2|\Delta u|^2\leq C(n)C^2_*+\frac{C(n)}{R^2}TC^2_*
\]
and therefore
\[
t^2|\Delta u|^2\leq C(n,K)A^2\left[1+T\left(1+\frac{1}{R^2}\right)\right]
\cdot\left(1+\frac{T}{R^2}\right)
\]
for any $(x,t)\in B_p(R)\times[0,T]$,
which completes the proof of the theorem.
\end{proof}


\section{Proof of Theorem \ref{the1}}\label{sect3}
In this section by using gradient and Laplacian estimates of
the previous section, we apply a maximum principle on complete
noncompact manifolds due originally to Karp and Li \cite{KaLi}
(see also Ni-Tam \cite{[Ni-Tam]}), to finish the proof of
Theorem \ref{the1}.

\begin{theorem}[Karp-Li \cite{KaLi} and Ni-Tam \cite{[Ni-Tam]}]\label{noncmax}
Let $(M,g)$ be an $n$-dimensional complete Riemannian manifold. Suppose
$f(x,t)$ is a smooth function on $M\times[0,T]$, $0<T<\infty$, such that
\[
\left(\frac{\partial}{\partial t}-\Delta\right)f(x,t)\leq 0
\quad\mathrm{whenever}\quad f(x,t)\leq 0.
\]
Let $f_+(x,t):=\max\{f(x,t),0\}$.  Assume that
\begin{equation}
\int^T_0\int_Me^{-ar^2(x)}f^2_+(x,t)d\mu dt\leq0
\end{equation}
for some constant $a>0$, where $r(x)$ is the distance to $x$ from
some fixed $p\in M$. If $f(x,0)\leq 0$ for all $x\in M$, then
$f(x,t)\leq 0$ for all $(x,t)\in M\times[0,T]$.
\end{theorem}

\begin{proof}[Proof of Theorem \ref{the1}]
Following the Hamilton's proof \cite{Hamilton} with a little modification.
We define $u_\epsilon=u+\epsilon$ satisfying
$\epsilon<u_\epsilon<A+\epsilon$ and the function
\begin{equation*}
\begin{aligned}
P(x,t)&:=t\left(\Delta u_\epsilon+\frac{|\nabla u_\epsilon|^2}{u_\epsilon}
\right)-u_\epsilon\left(n+4\log\frac{A}{u_\epsilon}\right).
\end{aligned}
\end{equation*}
By \cite{Hamilton}, we have
\begin{equation*}
\begin{aligned}
\left(\frac{\partial}{\partial t}-\Delta\right)P&\leq
-\frac{2t}{nu_\epsilon}\left(\Delta u_\epsilon-\frac{|\nabla u_\epsilon|^2}{u_\epsilon}
\right)^2+\left(\Delta u_\epsilon-\frac{|\nabla u_\epsilon|^2}{u_\epsilon}
\right)-2\frac{|\nabla u_\epsilon|^2}{u_\epsilon}.
\end{aligned}
\end{equation*}
By the Hamilton's arguments, we easily have that
\[
\left(\frac{\partial}{\partial t}-\Delta\right)P\leq 0
\quad\mathrm{whenever}\quad P\geq 0.
\]
This fact can be obtained by the following three cases.
\begin{enumerate}
\item If $\Delta u_\epsilon\leq\frac{|\nabla u_\epsilon|^2}{u_\epsilon}$, then we are done.
\item If $\frac{|\nabla u_\epsilon|^2}{u_\epsilon}\leq\Delta u_\epsilon\leq3\frac{|\nabla u_\epsilon|^2}{u_\epsilon}$, then we are also done.
\item If $3\frac{|\nabla u_\epsilon|^2}{u_\epsilon}\leq\Delta u_\epsilon$, then when
$P\geq 0$, we have
\[
2\left(\Delta u_\epsilon-\frac{|\nabla u_\epsilon|^2}{u_\epsilon}
\right)\geq\Delta u_\epsilon+\frac{|\nabla u_\epsilon|^2}{u_\epsilon}
\geq \frac{nu_\epsilon}{t}
\]
and hence we are done completely.
\end{enumerate}
Now, obviously, we have
\[
P(x,0)<0.
\]
By our assumptions on $u_\epsilon$, we also have
\[
P_+(x,t)\leq t\left(\Delta u_\epsilon
+\frac{1}{\epsilon} |\nabla u_\epsilon|^2\right),
\]
where $P_+(x,t):=\max\{P(x,t),0\}$. Thus, using
estimates \eqref{td1} and \eqref{td3gen},
for any $p\in M^n$, and $T, R>0$, we have
\begin{equation*}
\begin{aligned}
\int^T_0\int_{B_p(R)}&e^{-r^2(x)}P^2_+(x,t)d\mu dt\\
&\leq\int^T_0\int_{B_p(R)}e^{-r^2(x)}\left[t
\left(\Delta u_\epsilon+\frac{1}{\epsilon} |\nabla u_\epsilon|^2\right)\right]^2d\mu dt\\
&\leq\left(C(n)A+\frac{C(n)A^2}{\epsilon}\right)^2\int^T_0\int_Me^{-r^2(x)}
d\mu dt.
\end{aligned}
\end{equation*}
Since $Ric\geq 0$, by the Bishop
volume comparison theorem, we have that
\[
\int^T_0\int_{B_p(R)}e^{-r^2(x)}
d\mu dt<\infty.
\]
Then by letting $R\to\infty$, we conclude that
\[
\int^T_0\int_Me^{-r^2(x)}d\mu dt<\infty
\]
and hence
\[
\int^T_0\int_{B_p(R)}e^{-r^2(x)}P^2_+(x,t)d\mu dt<\infty.
\]
By the maximum principle for the complete noncompact manifold,
we conclude that $P(x,t)\leq 0$ for all $t\leq T$ and
hence the conclusion of Theorem \ref{the1} follows.
\end{proof}

\section{Proof of Theorem \ref{Kethe1}}\label{sect4}
The proof of Theorem \ref{Kethe1} follows from that of
Theorem 1 in \cite{Kotsch} with little modification,
but is included for completeness.
\begin{proof}[Proof of Theorem \ref{Kethe1}]
Let $H(x,y,t)$ be the heat kernel of the heat equation on $(M,g)$.
For any $t>0$ and $y\in M$, we set $u(x,s):=H(x,y,s+t/2)$, and
then $u$ is a smooth, positive solution to the heat equation on
$[0,T)$. By \cite{[Li-Yau]}, for any $\delta>0$, there exists a
constant $C_1=C_1(\delta)>0$ such that
\begin{equation}\label{keres7}
\frac{\exp\left(\frac{-d^2(x,y)}{(4-\delta)(s+t/2)}\right)}{C_1\mathrm{Vol}(B_y(\sqrt{s+t/2}))}
\leq u(x,s)\leq\frac{C_1}{\mathrm{Vol}(B_y(\sqrt{s+t/2}))}
\end{equation}
for all $x,y\in M$, and $s\geq 0$.

Letting
\[
A:=\frac{C_1}{\mathrm{Vol}(B_y(\sqrt{t/2}))},
\]
then the latter part of inequality \eqref{keres7} implies $u\leq A$ for all $x$ and
$s$. Since $Ric\geq0$, there exists a positive constant
$C_2:=C_2(n)$ such that
\[
\mathrm{Vol}(B_y(\sqrt{s+t/2}))\leq\mathrm{Vol}(B_y(\sqrt{t}))
\leq C_2\mathrm{Vol}(B_y(\sqrt{t/2}))
\]
for all $0\leq s\leq t/2$. Thus, by the front part of inequality \eqref{keres7}
and Theorem \ref{the1}, we have
\[
s\frac{\Delta u}{u}
\leq n+4\log\left(\frac{A}{u}\right)
\leq n+4\log(C^2_1C_2)+\frac{4d^2(x,y)}{(4-\delta)(s+t/2)}
\]
on $M\times [0,t/2]$. Setting $C=n+4\log(C^2_1C_2)$ and choosing at $s=t/2$,
from above, we conclude that
\[
(t/2)\frac{\Delta H}{H}(x,y,t)=(t/2)\frac{\Delta u}{u}(x,t/2)
\leq C+4\frac{d^2(x,y)}{(4-\delta)t}
\]
for all $x,y\in M$ and $t>0$.
\end{proof}


\section*{Acknowledgment}
The author would like to express his gratitude to the referee for careful readings
and many valuable suggestions.

\bibliographystyle{amsplain}

\begin{thebibliography}{30}
\bibitem{[Cala]}E. Calabi, An extension of E. Hopf's maximum principle with
an application to Riemannian geometry, Duke Math. J. 25 (1957), 45-56.

\bibitem{[CCG]}B. Chow, S.-C. Chu, D. Glickenstein,  C. Guentheretc,
J. Isenberg, T. Ivey, D. Knopf, P. Lu, F. Luo, L. Ni,  The Ricci flow:
techniques and applications. Part II: Analytic Aspects. Mathematical
Surveys and Monographs, AMS, Providence, RI, 2007.

\bibitem{[Gri]}A. Grigor'yan, Upper bounds of derivatives of the heat kernel
on an arbitrary complete manifold, J. Funct. Anal. 127 (1995), 363-389.

\bibitem{Hamilton}R.S. Hamilton, A matrix Harnack estimate for the heat equation.
Comm. Anal. Geom. 1 (1993), no.1, 113-126.

\bibitem{KaLi}L. Karp, P. Li, The heat equation on complete Riemannian manifolds,
http://math.uci.edu/~pli/, preprint, 1982.

\bibitem{Kotsch}B. Kotschwar, Hamilton's gradient estimate for the weighted kernel
on complete manifolds, Proc. AMS. 135 (2007), no.9, 3013-3019.

\bibitem{[Li-Yau]} P. Li, S.-T. Yau, On the parabolic kernel of the Schrodinger
operator, Acta Math. 156 (1986) 153-201.

\bibitem{[Ni1]}L. Ni, The entropy formula for linear heat equation, J. Geom. Anal.
14 (2004) 85-98. Addenda, 14 (2004), 369-374.

\bibitem{[Ni-Tam]} L. Ni, L.-F. Tam, K\"ahler-Ricci flow and the Poincar\'e-Lelong
equation. Comm. Anal. Geom. 12 (2004), no. 1-2, 111-141.

\bibitem{[Perelman]} G. Perelman, The entropy formula for the Ricci flow and
its geometric applications, arXiv:math.DG/0211159v1, 2002.

\bibitem{[Shi]}W.-X. Shi, Deforming the metric on complete Riemannian manifolds,
J. Differential Geom. 30 (1989), no. 1, 223-301.

\bibitem{Wu10} J.-Y. Wu, Li-Yau type estimates for a nonlinear parabolic equation on
complete manifolds, J. Math. Anal. Appl. 369 (2010), 400-407.

\bibitem{Yu}C.-J. Yu, Some estimates of fundamental solutions on noncompact mani-
folds with time-dependent metrics. To appear in Manuscripta Math., DOI:
10.1007/s00229-011-0517-y.

\end{thebibliography}

\end{document}